\documentclass[11pt]{article}
\usepackage{algorithm}
\usepackage{algpseudocode}
\usepackage{array}
\usepackage{makecell}
\usepackage{ragged2e}
\usepackage{lscape}	
\usepackage{graphicx}
\usepackage{appendix}
\usepackage{booktabs} 
\usepackage{tabularx}
\usepackage{comment}
\usepackage{subcaption}
\usepackage{amsmath, amsthm}
\usepackage{enumitem}

\newtheorem{Proposition}{Proposition}
    
\newcommand{\blind}{0}
	
    \makeatletter
    \renewcommand\section{\@startsection {section}{1}{\z@}%
                                       {-3.5ex \@plus -1ex \@minus -.2ex}%
                                       {2.3ex \@plus.2ex}%
                                       {\normalfont\fontfamily{phv}\fontsize{16}{19}\bfseries}}
    \renewcommand\subsection{\@startsection{subsection}{2}{\z@}%
                                         {-3.25ex\@plus -1ex \@minus -.2ex}%
                                         {1.5ex \@plus .2ex}%
                                         {\normalfont\fontfamily{phv}\fontsize{14}{17}\bfseries}}
    \renewcommand\subsubsection{\@startsection{subsubsection}{3}{\z@}%
                                        {-3.25ex\@plus -1ex \@minus -.2ex}%
                                         {1.5ex \@plus .2ex}%
                                         {\normalfont\normalsize\fontfamily{phv}\fontsize{14}{17}\selectfont}}
    \makeatother
\usepackage{multirow}

\usepackage{xargs}                      
\usepackage[pdftex,dvipsnames]{xcolor}  
\usepackage[colorinlistoftodos,prependcaption,textsize=tiny]{todonotes}
\newcommandx{\unsure}[2][1=]{\todo[linecolor=blue,backgroundcolor=blue!25,bordercolor=blue,#1]{#2}}
\newcommandx{\change}[2][1=]{\todo[linecolor=blue,backgroundcolor=blue!25,bordercolor=blue,#1]{#2}}
\newcommandx{\info}[2][1=]{\todo[linecolor=OliveGreen,backgroundcolor=OliveGreen!25,bordercolor=OliveGreen,#1]{#2}}
\newcommandx{\improvement}[2][1=]{\todo[linecolor=Plum,backgroundcolor=Plum!25,bordercolor=Plum,#1]{#2}}
\newcommandx{\thiswillnotshow}[2][1=]{\todo[disable,#1]{#2}}

	\usepackage{amsmath}
	\usepackage{graphicx}
	\usepackage{natbib} 
	\usepackage{url} 
    \usepackage{xcolor}
	
\begin{document}
		
		\def\spacingset#1{\renewcommand{\baselinestretch}%
			{#1}\small\normalsize} \spacingset{1}
		
		\if0\blind
		{
			\title{\bf Sensor-Driven Predictive Vehicle Maintenance and Routing Problem with Time Windows}
            
			\author{Iman Kazemian$^\dagger$ 
,   Bahar \c{C}avdar$^\ddagger$
            , and Murat Yildirim$^\dagger$ \thanks{Corresponding author: murat@wayne.edu}\\
			$^\dagger$ Industrial and Systems Engineering, Wayne State University, Detroit, US \\
             $^\ddagger$ Industrial and Systems Engineering, Rensselaer Polytechnic Institute, New York, US }
			\date{}
			\maketitle
		} \fi
		
		\if1\blind
		{

            \title
            {\bf Sensor-Driven Predictive Vehicle Maintenance
            and Routing Problem with Time Windows}
            
			\author{Author information is purposely removed for double-blind review}
			
\bigskip
			\bigskip
			\bigskip
			\begin{center}
				{\LARGE\bf{Sensor-Driven Predictive Vehicle Maintenance and \\ \vspace{0.3cm}Routing Problem with Time Windows}}
			\end{center}
			\medskip
		} \fi
		\bigskip
		
\vspace{-1cm}
	\begin{abstract}
 
Advancements in sensor technology offer significant insights into vehicle conditions, unlocking new venues to enhance fleet operations. While current vehicle health management models provide accurate predictions of vehicle failures, they often fail to integrate these forecasts into operational decision-making, limiting their practical impact. This paper addresses this gap by incorporating sensor-driven failure predictions into a single-vehicle routing problem with time windows. A maintenance cost function is introduced to balance two critical trade-offs: premature maintenance, which leads to underutilization of remaining useful life, and delayed maintenance, which increases the likelihood of breakdowns. Routing problems with time windows are inherently challenging, and integrating maintenance considerations adds significantly to its computational complexity. To address this, we develop a new solution method, called the Iterative Alignment Method (IAM), building on the structural properties of the problem. IAM generates high-quality solutions even in large-size instances where Gurobi cannot find any solutions. Moreover, compared to the traditional periodic maintenance strategy, our sensor-driven approach to maintenance decisions shows improvements in operational and maintenance costs as well as in overall vehicle reliability.

	\end{abstract}
			
	\noindent%



\section{Introduction} \label{s:intro}

In modern fleet management systems, sensor data is playing an increasingly critical role. The conventional role of sensor data has been to provide increased situational awareness for planners. To this end, sensor data has been used to generate real-time predictions on the current vehicle conditions and future failure risks - an area of research called vehicle health management. While these predictions have become indispensable over time, their contributions to vehicle operations management remained limited due to the difficulties associated with interpreting these predictions into operational decisions. A fundamental operational problem in fleet management is the Traveling Salesman Problem with Time Windows (TSPTW), which involves optimizing a single-vehicle route to minimize costs while ensuring that the vehicle visits customers within specified time frames. When sensor data and vehicle failure risk predictions are incorporated into the routing decisions, this problem evolves into the coordination of maintenance and operations, integrating TSPTW with predictive maintenance requirements. In this paper, we focus on this integration to make simultaneous routing and sensor-driven maintenance decisions.

In traditional applications of operations and maintenance, predictive maintenance schedules are planned based on predetermined intervals - \textit{think of the clich\'{e}, change oil once per 6,000 miles or 6 months of operation}. While straightforward to implement, this method relies on population-based reliability estimates, which are derived from general failure statistics for a fleet of vehicles. This approach fails to account for the unique operational conditions and specific wear patterns of individual vehicles. In practice, even identical vehicles often exhibit significant variability in degradation and failure risk trajectories due to factors such as manufacturing variability, material imperfections, driver behavior, and operational conditions. Although it is not feasible to directly observe these degradation processes, vehicle health management models utilize sensor data to infer them by developing fault signatures known as degradation signals. These signals correlate with the degradation severity, and they are used to predict current and future vehicle conditions and failure risks. Compared to reliability-based estimates, sensor-driven vehicle health management models provide significant insights into accurate failure risks and can be utilized to reduce failure risks and unnecessary maintenance.

Despite advancements in vehicle health management, conventional practices for routing and maintenance scheduling still predominantly rely on historical data and manual oversight, which lacks real-time adaptability and frequently leads to sub-optimal solutions (\cite{gackowiec2019general}). Although the value of sensor data and vehicle health management is recognized, their application has been largely confined to risk prediction, without extending these benefits to complex operational decisions. It remains an open challenge to see how these benefits can permeate into complex operational decisions. Addressing this challenge necessitates the development of a new generation of fleet operations models. These models should be capable of quantifying vehicle failure risks as predicted by vehicle health management systems and integrating these risk assessments with routing decisions, associated with the TSPTW. This integration aims to establish an end-to-end system that combines sensor-driven predictions with operational and maintenance decisions, thereby optimizing vehicle management processes.

In this work, we propose a framework for sensor-driven predictive vehicle maintenance and routing problem with time windows. The proposed method offers an integration of sensor-driven vehicle failure risks into decision optimization models for operations and maintenance. To ensure this integration, we use a cost quantification method through a dynamic maintenance cost function that translates failure risks to long-run average maintenance costs. This function establishes a trade off between (i) premature/early maintenance that underutilizes the equipment lifetime and (ii) postponed maintenance that increases failure risks. A unique feature of the dynamic maintenance cost function is its linkage to sensor-driven failure predictions. The sensor observation updates the failure risks and the associated dynamic cost function. This updated cost function is then integrated into the TSPTW model. By embedding the dynamic maintenance cost function into the TSPTW framework, we make simultaneous routing and vehicle maintenance decisions, leading to enhanced fleet efficiency. Our contributions can be summarized as follows:
\begin{enumerate}  
    \item[(i)] We present a new framework that integrates sensor-driven vehicle maintenance decisions with routing decisions for traveling salesman problem with time windows. This framework ensures that the benefits of vehicle health management models extend beyond mere predictions to enhance the efficacy of operational decisions, thereby reducing risks of vehicle failure and operational interruptions.
    \item[(ii)] We quantify sensor-driven failure risks through a dynamic maintenance cost function and embed it within the proposed operations and maintenance optimization model. This function uses sensor-driven predictions of vehicle failure risks to evaluate vehicle-specific long-run average maintenance costs; forming a link that enables strategic scheduling of operations and maintenance decisions. 
    \item[(iii)] To address the computational challenges associated with the integrated maintenance and routing problem, we develop a new algorithm, called the Iterative Alignment Method. The proposed method uses the structural properties of the integrated model to reformulate the problem into a series of routing problems.
    \item[(iv)] In addition, we conduct computational experiments to evaluate the performance of our solution method in different operational settings and compare its performance with traditional fleet management methods. In these experiments, we evaluate key performance metrics such as operational costs and the number of unexpected failures.
\end{enumerate}

The rest of this paper is organized as follows. In Section \ref{literature}, we review the related literature. In Section \ref{Sec:Description}, we present the problem description and our model along with structural results. We develop our solution method in Section \ref{Sec:Solution}. In Section \ref{experiments}, we present our computational experiments. Finally, in Section \ref{summary}, we conclude the paper.

\section{Literature Review}
\label{literature}

This literature review lays the groundwork for our study, which combines routing and maintenance. We first review studies on the Traveling Salesman Problem (TSP), along with its variants and applications. Then, we review maintenance scheduling, highlighting recent trends and new methods. Finally, we discuss studies that link routing and maintenance, providing context for our proposed approach.

TSP is a fundamental combinatorial optimization problem where the objective is to find a minimum-cost route that visits a set of nodes exactly once, returning to the starting point. TSP has several important variants applied across various fields, such as: 
the TSP with time windows (\cite{da2010general}), the Steiner TSP (\cite{rodriguez2019steiner}), the selective TSP (\cite{laporte1990selective}), the multi-objective TSP (\cite{psychas2015hybrid}), the multiple TSP (\cite{cheikhrouhou2021comprehensive}), and the bottleneck TSP (\cite{garfinkel1978bottleneck}). For further details on the current state of the art for these variants, please see the literature reviews on TSP and its variants (\cite{bock2024survey}; \cite{pop2024comprehensive}; \cite{toaza2023review}; \cite{pillac2013parallel}).

In recent years, there has been a significant emphasis on integrating TSP and its variants with other operational problems such as inventory management, scheduling, and technician routing and scheduling.  
For instance, combining TSP with inventory management allows for synchronized routing and inventory replenishment, which minimizes both travel distance and stockout risks (\cite{qiu2019optimal}). The integration of TSP with scheduling problems enhances the coordination between job allocation and routing, thus reducing total operational costs (\cite{he2023asymmetric}). Similarly, \cite{zamorano2017branch} addressed the multi-period technician routing and scheduling problem (MPTRSP) involving an external maintenance provider. The authors proposed a mixed integer programming model and a branch-and-price algorithm to solve this problem. \cite{irawan2017optimisation} proposed an optimization model for maintenance team routing and scheduling at offshore wind farms. \cite{si2022service} also considered a technician routing problem. They examined a leasing manufacturing system where the lessor manages the maintenance of the machinery. Their goal is to determine a route for the technicians which leads to minimum maintenance costs.

A significant operational problem that can be integrated with TSP is maintenance scheduling. Traditional maintenance strategies rely on fixed schedules or reactive repairs after failures \citep{coit2019evolution}. While these methods offer a simple operational strategy, they often overlook changing equipment conditions and operational environments, leading to unnecessary maintenance or unexpected breakdowns, increasing costs and downtime. To address this, Vehicle Health Management (VHM) \citep{benedettini2009state} has been used extensively by fleet operators to harness sensor data to monitor components in real-time and to predict failures. VHM provides significant visibility into the conditions of the vehicles and maintenance needs. There is a significant body of literature that focuses on integrating sensor-driven condition assessment into the maintenance scheduling problem, showcasing significant operational benefits \citep{shi2023stochastic,basciftci2020data,coit2019evolution}. However, this integration is missing for the joint maintenance and routing problem.

Due to inherent modeling and computational complexities, the existing approaches in routing and maintenance problems typically optimize these decisions sequentially. Typically, these approaches use a two-step process where the maintenance decisions are determined first, and routing is optimized around a set of fixed maintenance decisions \citep{lopez2016combined,fontecha2020combined}. 
Few studies optimized both decisions together. For example, Schindler integrated periodic maintenance into TSP using a geographic information system-based system, optimizing technician routes with algorithms to enhance planning and efficiency \cite{blakeley2003optimizing}.\cite{dhahri2015variable} developed a model for cases where vehicles undergo scheduled periodic maintenance at fixed times during their routes. They proposed a variable neighborhood search algorithm to minimize vehicle count and travel distance deviations. Similarly, \cite{dhahri2016vns} developed a mathematical model to determine whether maintenance should be performed at a node, focusing on identifying feasible maintenance times within the service time windows at each node. However, their model did not incorporate maintenance costs in its formulation.  Meanwhile, \cite{jbili2018integrated} proposed a vehicle routing and maintenance strategy for transcontinental transport, incorporating a random variable for failures and a fixed maintenance cost to minimize total costs of maintenance and late arrivals.

Despite significant visibility into vehicle conditions gained through the sensor-driven VHM applications, research combining routing and maintenance focused predominantly on periodic maintenance schedules that ignore the VHM inputs. In reality, insights gained from VHM can significantly improve operational outcomes, including routing and maintenance costs and operational reliability. This paper aims to bridge this gap by integrating routing and maintenance decisions using a sensor-driven maintenance approach.


\section{Problem Description}
\label{Sec:Description}

\begin{figure}[htbp]
    \centering\includegraphics[width=1\linewidth]{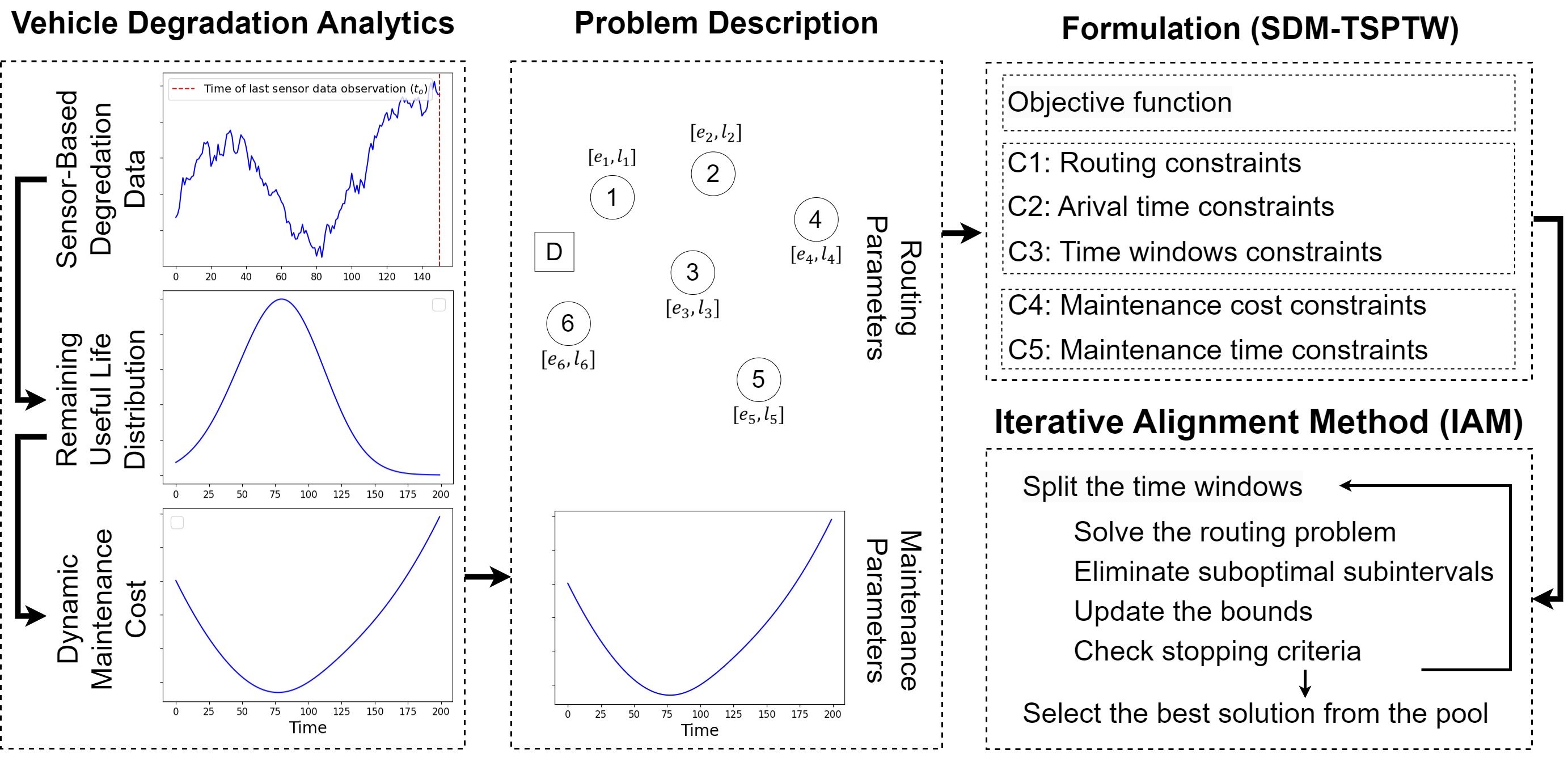} 
    \caption{Framework for integrated sensor-driven vehicle maintenance and routing problem.}
    \label{fig:ProblemOVerview}
\end{figure}

We consider a company providing a long-haul delivery service using a company-owned vehicle and performing the vehicle maintenance internally. 
Therefore, the company is making both operations and maintenance decisions.
On the operational side, the focus is on determining the optimal routes for the company-owned vehicle to satisfy time-sensitive demand from multiple customers. On the maintenance side, the focus shifts to the scheduling of maintenance at maintenance-capable locations to reduce failure risks and minimize maintenance costs. In this setting, routing and maintenance decisions are highly interdependent. For example, improper maintenance decisions may lead to an unexpected vehicle breakdown; on the other hand, the routing decisions affect the vehicle degradation level and the maintenance decisions as a result. 
In this paper, we study how to make integrated routing and maintenance decisions, considering the dynamic relationship between them. We present the details of the underlying problem environment below. Figure \ref{fig:ProblemOVerview} displays an overview of the problem elements and solution framework.

In our setting, we consider a company that provides a specialized service to its customers. The service is performed by a traveling crew.
We assume that customer requests are finalized before the service vehicle is dispatched from the depot; therefore, the demand is deterministic for the routing decisions. Let $n$ be the number of customer locations to be visited and $\{1, 2, \dots, n\}$ represent the set of customer nodes.
We represent the road network covering the customer locations and the depot by $G=(N,A)$, where $N$ is the set of nodes including the depot (node 0) and customers, and $A$ is the set of edges between each pair of nodes. 
We denote the minimum travel time between nodes $i$ and $j$ by $d_{ij}$. Each customer order is associated with a hard time window $[e_i, l_i]$, where $e_i$ is the start and $l_i$ is the end of the time window.
Therefore, the routing decisions resemble the TSPTW.

In addition to the routing decisions, we also consider the degradation of the vehicle during the route. 
We assume that the vehicle deteriorates as it travels a long-haul route, which may necessitate a maintenance operation to ensure operational continuity.
If the maintenance need arises en route, the maintenance operation can be performed by the travel crew or a third-party team. Regular maintenance operations are called \textit{preventive maintenance}. 
We assume that the preventive maintenance duration, denoted by $p$, is deterministic, and the vehicle can continue with the planned route after the maintenance with no further delay.
However, if the vehicle breaks down unexpectedly, then the vehicle needs to go through a \emph{corrective maintenance}, which is much more costly.
We assume that a subset of customer locations are suitable for performing vehicle maintenance. 
We call them the \textit{maintenance-capable} nodes and denote the set by $N^{'}$, where $N^{'}\subseteq N$. The nodes in $N^{'}$ are distinct from $N\setminus N^{'}$ either due to providing the necessary infrastructure for vehicle maintenance or through agreements with that specific customer or maintenance provider. 
Table \ref{Tab:Notation} provides a summary of our notation.

Maintenance scheduling decisions alone inherent significant tradeoffs. 
If maintenance is performed too early, then the company cannot make efficient use of the vehicle's useful life and incurs unnecessarily high maintenance costs.
On the other hand, if maintenance is performed too late, then unexpected vehicle breakdowns may occur, resulting in high corrective maintenance costs and other additional costs due to disrupted operations, e.g., paying penalties for unmet customer orders.
Therefore, maintenance decisions need to consider both aspects.
Our goal is to make joint routing and maintenance decisions to minimize the total cost while serving customer orders with time windows considering the vehicle sensor data on degradation level. We call this problem the Sensor Driven Maintenance and Traveling Salesman Problem with Time Windows (SDM-TSPTW). 
To integrate these two decisions, we model the relationship between the routing and maintenance decisions. More specifically, we determine the degradation level of the vehicle based on its use, which is based on the travel duration determined by the route. For this purpose, we first derive the dynamic maintenance cost function in the following subsection.

\begin{table}[h!]
\centering
\caption{Summary of Notation.}
\begin{tabularx}{\textwidth}{>{\raggedright\arraybackslash}p{0.4\linewidth}>{\raggedright\arraybackslash}p{0.55\linewidth}}
    \toprule
    \textbf{Notation} & \textbf{Description} \\
    \midrule\( N=\{0,1,\dots,n\} \) & Set of nodes ($n$ customers and a depot) \\
    \( N^{'} \in N \) & Set of maintenance nodes \\
    \( d_{ij} \) & Travel time between nodes \(i\) and \(j\) \\
    \( [e_i, l_i] \) & Time window for node \(i\) \\
    \( p \) & Duration of planned maintenance \\
    \( C^r \) & Cost of unit travel time \\
    \( C^p \) & Cost of planned maintenance \\
    \( C^f \) & Cost of failure \\
    \( t_o \) & Time of the sensor data observation \\
    \( f_{t_o}(t) \) & Dynamic vehicle maintenance cost function after duration \(t\) since the sensory data observation \(t_o\) \\
    \( T_{min} \) & Time at which the dynamic maintenance cost is at its minimum, \(T_{min}=\underset{t}{\text{argmin}} f_{t_o}(t)\) \\
    \( \lambda \) & Minimum dynamic maintenance cost, \( \lambda=\underset{t}{\text{min}} f_{t_o}(t) \) \\
    \bottomrule
\end{tabularx}
\label{Tab:Notation}
\end{table}


\subsection{Linking Vehicle Degradation with Routing Decisions}
\label{Sec:Link}
An important aspect of our modeling framework is the integration of sensor-driven degradation analytics into routing and maintenance decisions. To enable this integration, we develop a framework that achieves two linked goals: (i) leveraging sensor-driven insights to develop accurate predictions on vehicle degradation and failure risks and (ii) translating these risks into maintenance costs considering the trade-offs. In line with these goals, we provide a two-stage process, each corresponding to one of these goals.

\textbf{Stage 1: Vehicle Degradation Modeling and Remaining Life Prediction:} 

The first stage involves a prediction task that uses the sensor data from the vehicle to predict the remaining life distribution. Unlike classical reliability-based formulations that focus on population-wide failure statistics, our proposed failure prediction model leverages on degradation insights that can be acquired from the sensor data. Specifically, we are looking for indicators within the sensor observations that correlate with evolving degradation processes. These sensor data sources could be physical indicators such as vibration, which increases with issues like wear/tear and misalignment; or electrical indicators such as voltage and current readings which exhibit significant trends with certain degradation processes (e.g. battery). The correlations between the sensor data and the degradation level are often not explicit. To derive this correlation, condition monitoring uses multi-variate sensor-data from vehicles, and builds dimensionality reduction and predictive models (e.g. principal component analysis (PCA), Kalman filtering, Bayesian inference) to infer the severity of degradation, typically through a single unified measure called the degradation signal \citep{gebraeel2005residual}. Degradation signal constitutes a holistic estimate of the vehicle degradation, which can be used to make predictions on current state of the vehicle health and future risks of failure. The sensor signals also contain several sources of inevitable uncertainties, due to material imperfections, manufacturing variations, sensor errors, and inherent uncertainty of the degradation processes. Therefore, they are typically modeled as stochastic processes.

In our framework, $D(t)$ denotes the degradation signal of the vehicle at time $t$, which evolves as a continuous-time continuous-state stochastic process, as follows:
\begin{equation}
D(t) = \phi(t;\kappa,\theta) + \epsilon(t;\sigma),
\end{equation}
where $\phi$ is a parametric degradation function that is governed by two parameters. The term $\kappa$ is a deterministic parameter that is typically used to denote population-based degradation characteristics. The second term $\theta$ is a stochastic parameter that characterizes vehicle-specific degradation characteristics. We commonly do not have perfect information about parameter $\theta$. However, we can use the sensor data to improve our estimation accuracy for this parameter using a Bayesian learning approach. 
We assume that there is prior knowledge on the distribution of this parameter through past observations and engineering knowledge, which we  define as the prior distribution for this parameter, denoted by $\pi(\theta)$. After observing the degradation signal, we refine this distribution. Let $d_{t_o}$ be the set of degradation signal observations at the time of observation  $t_o$. Upon observing the signal $d_k$, we can compute the sensor-updated posterior distribution of  the stochastic parameter $\theta$, denoted by $u(\theta)$, as follows:
\begin{equation}
u(\theta) = \frac{P(d_{t_o} | \theta) \pi(\theta)}{P(d_{t_o})} \label{eq:U}.
\end{equation}

Once we obtain an improved estimate on the degradation parameter, we compute predictions on the remaining life distribution. To do so, we first define the failure condition. as the first time the degradation signal $D(t)$ reaches a failure threshold $\Lambda$. Using this failure condition definition, we can use the sensor-updated degradation signal parameters to predict the remaining life of a vehicle of age $t_o$, denoted by $R_{t_o}$, using the following steps: 
\begin{itemize}[noitemsep,nolistsep,leftmargin=4mm]
\item \textit{Step1:} Select a sufficiently large number $M$. 
\item \textit{Step2:} For each $i = 1, 2\dots, M$:
\begin{itemize}[noitemsep,nolistsep,leftmargin=5mm]
\item \textit{Step 2.1:} Generate a realization of parameter $\theta$, denoted by $\tilde{\theta}_{i}$, using its posterior distribution $u(\theta)$. 
\item \textit{Step 2.2:} Using the realized $\tilde{\theta}_{i}$, simulate $\{D(\tau+t_o|\tilde{\theta}_{i}, D(t_o)), \; \forall \tau \geq 0 \}$. In this simulation, we condition on two sources of information: (i) $\tilde{\theta}_{i}$, the realization of the stochastic parameter, and (ii) $D(t_o)$, the degradation signal amplitude at the time of observation $t_o$ (e.g., most recently observed degradation severity). 
\item \textit{Step 2.3:} Compute the $i^{th}$ realization for the remaining life denoted by $\tilde{r}_{t_o,i}$ by calculating the shortest time required for the simulated degradation signal realization to exceed the failure threshold, i.e., 
$\tilde{r}_{t_o,i} = min\{\tau>0| D(\tau+t_o|\tilde{\theta}_{i}, D(t_o))\geq \Lambda\}$. 
\end{itemize}
\item \textit{Step 3:} Collect realizations of remaining life across all samples $\tilde{r}_{t_o,i}, \; \forall i \in \{1,\cdots,M\}$ to estimate a distribution for the remaining life corresponding to the random variable $R_{t_o}$.
\end{itemize}

For certain special cases of the degradation models, there may also be closed-form approximations for the remaining life distribution (see \cite{gebraeel2005residual} for an example).

\textbf{Stage 2 - Dynamic Cost Calculation:} Given the remaining life distribution prediction from Stage 1, we derive a dynamic maintenance cost function. The cost function integrates the sensor-driven remaining life distribution predictions into a cost rate calculation that aims to balance the trade-off between (i) premature maintenance that underutilizes the equipment lifetime, and (ii) late maintenance that increases the risk of failure. We formulate the cost rate as a function of the vehicle maintenance time $t$, based on the sensor data observed at time $t_o$, $f_{t_o}(t)$, as follows:

\begin{equation}
f_{t_o}(t)=\frac{P(R_{t_o}>t) C^p+P(R_{t_o}<t) C^f}{\int_0^t P(R_{t_o} > z) \, dz + t_o} \label{eq:MaintenanceCost}.
\end{equation}
where, the terms $C^p$ and $C^f$ represent the cost of preventive maintenance and failure operations, respectively. The parameter $C_f$ represents the combined cost of corrective maintenance and the additional cost incurred due to operational disruptions. The numerator computes the total expected maintenance  cost per maintenance cycle. The denominator calculates the expected length of the maintenance cycle \citep{ elwany2008sensor}.

\subsection{Formulation}

Having developed the dynamic maintenance cost function, we now present our mathematical model for SDM-TSPTW to make joint routing and maintenance decision and minimize the total routing and maintenance-related costs. Our model integrates the dynamic maintenance cost function in Equation (\ref{eq:MaintenanceCost}) to account for the sensor-driven information regarding the health of the equipment. Our primary decision variables are as follows:
\begin{itemize}
\item  $x_{ij}$: binary variable that equals $1$ if node $j$ is visited immediately after node $i$, or $0$ otherwise,
\item $m_i$: binary variable that equals $1$ if maintenance is performed at node $i$, or $0$ otherwise.
\end{itemize}
Based on these variables, we also compute the following auxiliary variables:
\begin{itemize}

\item  $y$: binary variable that equals $0$ if maintenance is performed during the route, or $1$ otherwise,
\item $u_i$: arrival time at node $i$,
\item $\tau$: total travel and maintenance duration,
\item $\pi$: time of the maintenance, if it is performed during the route,
\item $\gamma$: cost of the maintenance, if it is performed during the route.
\end{itemize}
The mathematical formulation for the proposed SDM-TSPTW is as follows:
\small{
\begin{align}
\text{Min } z &= C^r \tau + \gamma + y \lambda \label{eq:obj} \\
\sum_{j=0}^n x_{ij} &= \sum_{j=0}^n x_{ji} & \forall i \in N \label{eq:continiuty} \\
\sum_{\substack{j=0 \\ j\neq i}}^n x_{ij} &= 1 & \forall i \in N \label{eq:mustleave} \\
u_0 &= 0 & \label{eq:uInitialize} \\
u_{j} &\geq u_i + d_{ij} + m_{i} p - M(1-x_{ij}) & \forall i \in N, j \in N-\{0\} \label{eq:arrivaltime} \\
l_{i} &\geq u_i \geq e_i & \forall i \in N-\{0\} \label{eq:late} \\
\tau &\geq u_i + d_{i0} + m_i p - M(1- x_{i0}) & \forall i \in N-\{0\} \label{eq:routedur} \\
\pi &\geq u_{i} - (1-m_{i})M & \forall i \in N^{'} \label{eq:mainTime1} \\
\pi &\leq u_{i} + (1-m_{i})M & \forall i \in N^{'} \label{eq:mainTime2} \\
\gamma &\geq f_{t_o}(\pi) - My &  \label{eq:maintCost} \\
y &\geq 1-\sum_{i \in N^{'}}m_i & \label{eq:noMaint} \\
M \sum_{i \in N^{'}}m_i &\geq \tau-T_{min} & \label{eq:mustMaint} \\
x_{i,j} &\in \{0,1\} & \forall i,j \in N \label{eq:typeStart}\\
m_{i} &\in \{0,1\} & \forall i \in N \\
u_{i} &\geq 0 & \forall i \in N \\
y &\in \{0,1\} & \\
\pi, \gamma, \tau &\geq 0 & \label{eq:typeEnd}
\end{align}
}
\normalsize

The objective function in \eqref{eq:obj} computes the total cost. The first term is the cost of time spent during the service delivery, including the maintenance duration if it is performed. The second and the third terms compute the cost of maintenance. If the maintenance is performed during the route, then cost $\gamma$ is incurred based on the dynamic maintenance cost function, which is computed in the constraints. Otherwise, we assume that the vehicle maintenance is performed at a future time (outside of the current route) at its minimum cost $\lambda$.
Constraints \eqref{eq:continiuty} ensure the continuity of the route. Constraints \eqref{eq:mustleave} enforce that each customer node is visited exactly once. Vehicle dispatching time is set to 0 in constraint \eqref{eq:uInitialize}.  
Arrival time at each node is computed in constraints \eqref{eq:arrivaltime}.
Constraints \eqref{eq:late} enforce the customers' time windows. 
Constraints \eqref{eq:mainTime1} and \eqref{eq:mainTime2}, where $M$ is a large value, compute the time of the maintenance $\pi$. Constraint 
\eqref{eq:maintCost} computes a lower bound for the maintenance cost at time $\pi$, if a maintenance is performed during the route (i.e, $y=0$). Since we are minimizing, $\gamma$ will be equal to its lower bound. Constraint \eqref{eq:noMaint} determines whether a maintenance is performed during the route or not.
Constraint \eqref{eq:mustMaint} enforces a maintenance if the route duration extends \( T_{min} \), the time at which the dynamic maintenance cost is at its minimum. This is by our assumption that if we reach the minimum maintenance cost during the travel, we perform the maintenance. Constraints \eqref{eq:typeStart}-\eqref{eq:typeEnd} specify the type of the decision variables.

Solving this problem using commercial solvers poses limitations both in terms of the computation time and size of the instances.  
Hence, our efforts are in developing efficient solution methods. For this purpose, we first develop some structural results to identify a set of suboptimal solutions. Our solution method in Section~\ref{Sec:Solution} is based on these results.

\subsection{Structural properties} 

Consider a maintenance node $i$ with time window $[e_i, l_i]$. One can split this time window into subintervals as shown in Figure~\ref{fig:first}. Let $q$ and $k$ be indices for the subintervals, i.e., subinterval $q$ and subinterval $k$, and subinterval $q$ denoted by $[h_{q-1},h_q]$.  
Accordingly, we define the following notation:
\begin{itemize}
    \item Let \(\tau^i_q\) be the optimal total route duration (including the travel and maintenance), when the time window of maintenance node $i$ is restricted to its subinterval \(q\).
    \item  Let \(\overline{g}^i_q\) and \(\underline{g}^i_q\) be the maximum and minimum values of the maintenance cost that can be achieved in time window subinterval \(q\) for maintenance node $i$, respectively. That is, $\overline{g}^i_q = max\{ f_{t_o}(s) | h_q \geq s \geq h_{q-1} \}$ and $\underline{g}^i_q = min\{ f_{t_o}(s) | h_q \geq s \geq h_{q-1} \}$.

\end{itemize}

For the following propositions, consider the set of disjoint time window subintervals for node $i$, denoted by $S_i$.

\begin{Proposition}
\label{Prop:dominance2}
If  \(C^r \tau^i_q + \underline{g}^i_q > C^r \tau^i_k + \overline{g}^i_k\), then the optimal solution to the SDM-TSPTW cannot reside within time window subinterval \(q\) of node $i$.
\end{Proposition}

\begin{Proposition}
\label{Prop:upper}
\(\min_{q \in S_i}\{C^r \tau^i_q+\overline{g}^i_q\}\) is an upper bound for the SDM-TSPTW  given that a maintenance is enforced to be performed at node $i$.
\end{Proposition}

\begin{Proposition}
\label{Prop:lower}
\(min_{q \in S_i}\{C^r \tau^i_q + \underline{g}^i_q\}\) is a lower bound for the SDM-TSPTW given that a maintenance is enforced to be performed at node $i$.
\end{Proposition}

The proofs are in Appendix \ref{appendix:proof}. Using Propositions \ref{Prop:dominance2}-\ref{Prop:lower}, we can calculate lower and upper bounds for the objective function value of the SDM-TSPTW for each maintenance node $i \in N^{'}$ considering subintervals of its time window and subsequently eliminate the dominated subintervals according to these bounds.  We use these properties in the development of our solution methodology, outlined in the next section.

\section{Solution Methodology}
\label{Sec:Solution}
 
The idea of our solution method lies in the iterative elimination of suboptimal time-window subintervals of maintenance nodes. In each iteration, we consider a all the maintenance nodes, $i\in N'$, and for each of time window subintervals, $q$, we solve the TSPTW problem assuming that node $i$ can be visited only during subinterval $q$, ignoring the maintenance cost. Then, using the corresponding routing cost and the minimum maintenance cost in this time interval, i.e., $\underline{g}^i_q$, we can determine whether time window subinterval $q$ for node $i$ is dominated or not under the current upper bound.
During the iterative process, we adjust the time windows of the maintenance-capable nodes progressively. In each iteration, we tighten the time window subintervals, by creating further subintervals, and integrate the corresponding dynamic maintenance cost into the objective function. Meanwhile, we determine the routing decisions using a popular solution method, called the Lin-Kernighan-Helsgaun (LKH) algorithm. Using this approach, we can reduce the complexity of the sensor-driven predictive maintenance decisions and solve the problem without initially considering the maintenance cost.  
We present the detailed solution methodology as follows:

\textbf{Step 1-Initialization:} We begin with initializing the problem input, including the network, time windows, and dynamic maintenance cost function $f_{t_o}(t)$, which is constructed based on the degradation signal.

\begin{figure}[h]
     \centering
     \begin{subfigure}[b]{0.49\textwidth}
         \centering
 \includegraphics[width=1\linewidth]{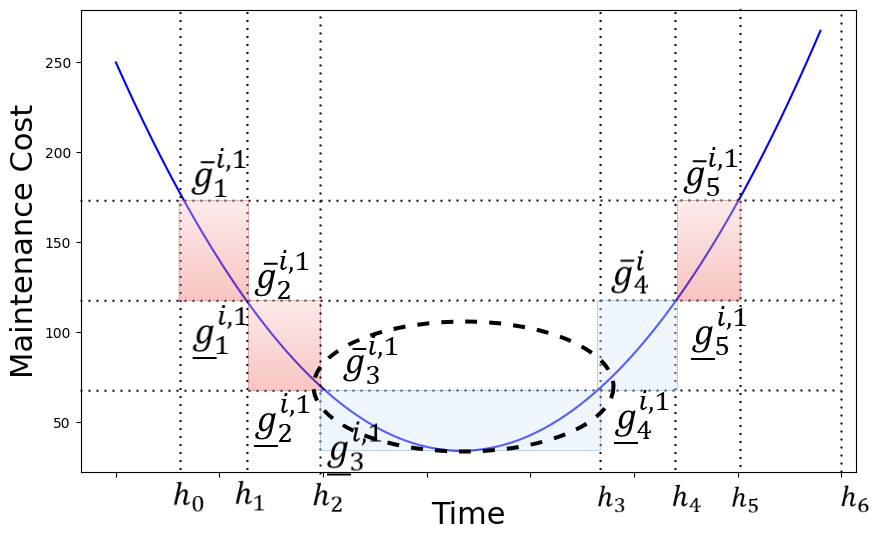}
         \caption{Time window subintervals for node $i$ based on its original time window $[e^1_i, l^1_i]$ created in iteration 1 for $b=5$. Note: $h_0=e^1_i$ and $h_5=l^1_i$.}
         \label{fig:first}
     \end{subfigure}
     \hfill
     \begin{subfigure}[b]{0.49\textwidth}
         \centering
         \includegraphics[width=1\linewidth]{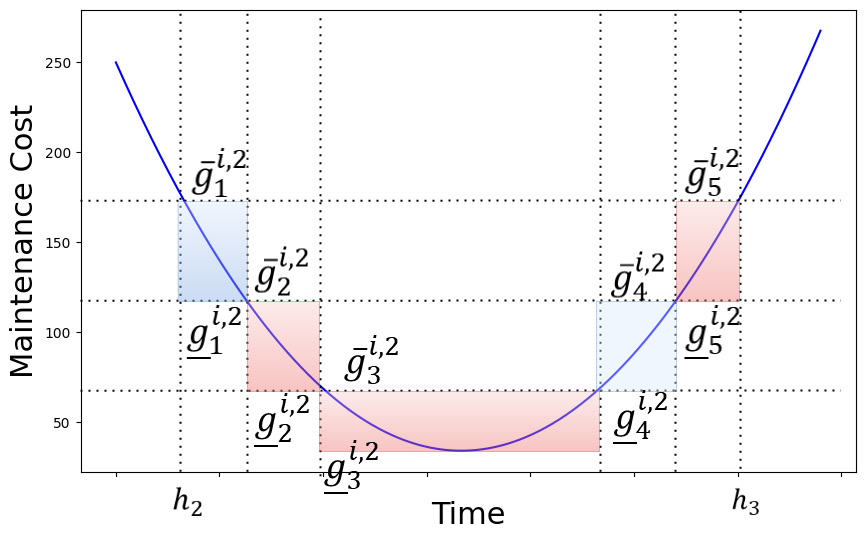}
         \caption{Recalculated subintervals for $[h_2, h_3]$ for node $i$ in iteration 2 for $b=5$.\\}
         \label{fig:recalc}
     \end{subfigure}
       \caption{\small{Dynamic maintenance cost function within time-window subintervals of an arbitrary node \(i\).  
    In subfigure (a), the time window \([e^1_i, l^1_i]\) is initially divided into five subintervals. Subintervals highlighted in red demonstrate the ones shown to be suboptimal, blue subintervals survive this iteration and are split further in the next iterations.
    Subfigure (b) zooms onto the subinterval time window $[h_2, h_3]$ in iteration 2 and displays its five subintervals in this iteration.
    }}
    \label{fig:DMC_combined}
\end{figure}

\textbf{Step 2-Iterative Routing Decisions with Dynamic Time Window Adjustment:} In this step, we iteratively create time-window subintervals for maintenance-capable nodes and refine solution by removing the suboptimal ones. The process is as follows:
\begin{enumerate}[label=(\roman*)]
    \item \textbf{Creating Time Window Subintervals:} For each maintenance-capable node \(i\), we split its current time window into subintervals. 
    Subintervals are created so that the differences between the maintenance costs at the start and the end of the time window subintervals are equal. 
    As a result, subintervals may vary in length, depending on how the maintenance cost function changes over the time window. However, the maintenance cost change over each subinterval is constant.

To create time window subintervals over which the changes in dynamic maintenance cost are the same, we first compute a constant $\delta$ as follows.
\begin{equation}
\delta =
\begin{cases}
\frac{|f_{t_o}(e_i^v) - f_{t_o}(l_i^v)|}{b}, \qquad \text{if $f_{t_o}(h)$ is monotone over $[e^v_i, l^v_i]$,}\\[0.5cm]
\frac{| f_{t_o}(e_i^v)- \underset{h\in [e^v_i, l^v_i]}{\min} f_{t_o}(h)| + |f_{t_o}(l_i^v)- \underset{h\in [e^v_i, l^v_i]}{\min} f_{t_o}(h)| }{b}, \qquad\text{otherwise.}\\
\end{cases}
\end{equation}

    Subsequently, time-window subintervals $[h_{q-1}, h_{q}]$, where $q\in \{1, 2, \dots, b\}$ are determined as follows: $h_0$ is set to be $e_i^v$, and starting with $q=1$, $h_q$ is chosen such that $f_{t_o}(h_q)-f_{t_o}(h_{q-1})= \delta$.  This process is illustrated in Figure~\ref{fig:first}, where we display the initial subintervals created by splitting the time window of node $i$ into $b=5$ subintervals, where $h_0=e^v_i$ and $h_5=l^v_i$.

    This method is preferred over choosing subintervals of similar length because it inherently incorporates the shape of the maintenance cost function \(f_{t_o}(h)\) into the process. When  \(f_{t_o}(h)\) is flatter, the maintenance cost is not as sensitive and larger intervals are created; when it is steeper, smaller intervals are used for better control. 

    \item \textbf{Adjusting the Distance Matrix:} When considering to perform maintenance at node \(i\), we temporarily update the travel times from node \(i\) to all other accessible nodes \(j\) by adding the maintenance duration $p$ to the original travel durations.

    \item \textbf{Solving the Adjusted TSPTW:} After adjusting the travel times, we use the LKH algorithm (\cite{helsgaun2009general}) to solve the TSPTW for the current time window subintervals of node $i$. Note that, currently, the maintenance is assumed to be performed at node $i$.
    \item \textbf{Eliminating Suboptimal Time Window Subintervals:} Once we obtain the routing solutions with LKH, we detect the suboptimal time window subintervals among all surviving ones using Propositions \ref{Prop:dominance2}-\ref{Prop:lower}, and we remove them from further considerations.  If maintenance is not performed during the route, the lower and upper bounds for the maintenance cost function uses the minimum dynamic maintenance cost, denoted by \(\lambda\).
    
    In Figure~\ref{fig:DMC_combined}, we present a demonstration for creating time window subintervals. Let us assume that time interval subintervals $q=1,2,5$ that are highlighted in red in Figure \ref{fig:first} are marked to be suboptimal and hence removed. Subintervals $q=3,4$ survive the elimination and will be reconsidered later. When we are processing node $i$ again, its remaining time window subintervals will be further siplit. For example, in Figure \ref{fig:recalc}, we display further splitting of the surviving time window subinterval $[h_2, h_3]$ into $b=5$ new intervals. 
    
    \item \textbf{Recalculating Bounds:} With each new routing solution, we update lower and upper bounds for the optimal objective function value using Propositions \ref{Prop:upper} and \ref{Prop:lower} as follows.
        \begin{itemize}
        \item \textbf{Upper bound - iteration $v$:} \(U^v\)= \(\min\{U^{v-1}, \min_{i \in N'} \{\min_{q \in S^{v-1}_i}\{C^r \tau^i_q + \overline{g}^i_q\}\}\}\),
        \item \textbf{Lower bound - iteration $v$:} \(L^v\)=  \(\max\{L^{v-1}, \min_{i \in N^{'}}\{\min_{q \in S^{v-1}_i}\{C^r \tau^i_q + \underline{g}^i_q\}\}\}\),
    \end{itemize}
    
    where $S^{v}_i$ is the set of time window subintervals for node $i$ in iteration v.

    Steps described in (i)-(v) are iterated over the maintenance-capable nodes, where a node can be iterated over several times with further time-window splitting, until the stopping criteria is met.

    \item \textbf{Termination:} We stop the iterations when the gap between upper and lower bounds is less than some small $\epsilon$. Once the stopping criterion is met, we choose the solution with the minimum total travel and maintenance cost.
\end{enumerate}

We call our solution method the \emph{Iterative Alignment Method} (IAM).

\begin{Proposition}
\label{Prop:optimality}
For any given $\epsilon$, the Iterative Alignment Method (IAM) finds solutions within $\epsilon$ of the optimal objective function value in finite number of iterations.

\end{Proposition}

 The proof is provided in Appendix~\ref{appendix:proof}.

\begin{algorithm}[h]
\caption{Iterative Alignment Method (IAM) for SDM-TSPTW}
\label{Alg:Heuristic}
\begin{footnotesize}
\begin{algorithmic}[1]
\Require Instance data and dynamic maintenance cost function.
\Ensure Integrated routing and maintenance decisions
\State Initialize the set of demand nodes \(N\), time windows \([e_i, l_i]\), and dynamic maintenance cost function \(f_{t_o}(h)\). Set the initial iteration number $v \leftarrow 0$ and initialize $\delta^0$.
\Repeat
    \For{each maintenance-capable node \(i\) in \(N'\)}
        \State Split the time window \([e^v_{i,o}, l^v_{i,o}]\) for every undominated subintervals $o \in S_i^v$ into subintervals \([h_{q-1}, h_q]\), $q \in \{1,2,...,b\}$, such that the difference in maintenance costs between the start and end times of each subinterval is constant.
        \State Calculate \(\delta^v\)= $
        \begin{cases}
\frac{|f_{t_o}(e^v_{i,o}) - f_{t_o}(l^v_{i,o})|}{b}, \qquad \text{if $f_{t_o}(h)$ is monotone over $[e^v_{i,o}, l^v_{i,o}]$,}\\[0.2cm]
\frac{| f_{t_o}(e^v_{i,o})- \underset{h\in [e^v_{i,o}, l^v_{i,o}]}{\min} f_{t_o}(h)| + |f_{t_o}(l^v_{i,o})- \underset{h\in [e^v_{i,o}, l^v_{i,o}]}{\min} f_{t_o}(h)| }{b}, \qquad\text{otherwise.}\\
\end{cases}
        $ .
            
        \State \State Determine time-window subintervals \([h_{q-1}, h_q]\) for \(q \in \{1, 2, \dots, b\}\) by setting \(h_0 = e^v_{i,o}\) and iteratively choosing \(h_q\), so that maintenance cost change is less than or equal to \(\delta^v\) for every subinterval.    
        \For{each subinterval \(q\) of node \(i\)}
            \State Adjust the arrival times from node \(i\) to all other accessible nodes \(j\) to account for maintenance at node \(i\) during subinterval \(q\).
            \State Solve the adjusted TSPTW using the LKH algorithm.
            \State If LKH return a solution store the solution and corresponding objective value \(\tau^i_q\).
        \EndFor
    \EndFor
    \For{each maintenance-capable node \(i\) in \(N'\)}
        \For{each pair of subintervals \(q\) and \(k\) for node \(i\)}
            \If{\(C^r \tau^i_q + \underline{g}^i_q \geq C^r \tau^i_k + \overline{g}^i_k\)}
                \State Eliminate subinterval \(q\) from further consideration.
            \EndIf
        \EndFor
    \EndFor
    \State Update the upper bound $U^v$ as \(\min\{U^{v-1},\min_{i \in N'} \{\min_{q \in S^{v-1}_i}\{C^r \tau^i_q + \overline{g}^i_q\}\}\}\).
    \State Update the lower bound $L^v$ as \(\max\{L^{v-1}, \min_{i \in N'}\{\min_{q \in S^{v-1}_i}\{C^r \tau^i_q + \underline{g}^i_q\}\}\}\).

    \State 
    Calculate the incurred total costs $z_{i,q}^{v-1}$ for each node and each subinterval that was not eliminated in line 17. Put these solutions into the current solution pool.
    \State
    Update \(\delta^{v+1} \leftarrow\delta^v/b\), $v \leftarrow v+1$
    
\Until{$U^v - L^v \leq \epsilon + 2 \delta^v$}
\State Check the solutions from the solution pool in iteration $v-1$ and choose the best solution in terms of minimum total cost $z$.
\end{algorithmic}
\end{footnotesize}
\end{algorithm}


\section{Computational Experiments}
\label{experiments}
In this section, we present a series of computational experiments to showcase the performance of our proposed SDM-TSPTW approach. We evaluate the performance of our approach in two parts. In the first part, we compare our IAM solution method with a benchmark method that uses Gurobi, and we evaluate our performance in terms of solution quality and computation time. In the second part, we shift our focus to the evaluation of operational outcomes. We compare the proposed sensor-driven approach with a traditional periodic maintenance-based benchmark model that does not use sensor information. We showcase how the additional insights gained from sensor-data can be used to attain better operational outcomes, as manifested through significant reductions in routing and maintenance costs.

In all our computational experiments, our problem instances are based on the TSP with time windows instances by \cite{gendreau1998generalized}.
These instances use the following naming convention: nXwY.Z, where X is the number of nodes, Y indicates the width of the time windows, and Z is the index of the instance. We use instances with up to 80 nodes. Instance names are reported on the tables where we present the numerical results.
We choose maintenance-capable nodes randomly using a p-median problem.
Our dynamic maintenance cost function is calculated based on a degradation experiment for rotating machinery. An accelerated life testing experiment was conducted for a population of rotating machinery, and the vibration data was continuously collected throughout the lifetime to construct the degradation signals. These degradation signals were used to predict remaining life distribution. For a more detailed explanation on the data and life prediction procedure, please see \cite{gebraeel2005residual}. The resulting remaining life predictions were used to evaluate dynamic maintenance cost function as shown in Section 3.1. 
To convert the travel time to cost, we use the conversion rate based on the guidelines from the U.S. Department of Transportation, where one unit of travel has a cost of 0.72 units \citep{bts_transportation_spending}. We implement the IAM algorithm using Python. We use Gurobi to solve the related optimization problem to obtain the lower bounds as described in the following subsection.
All computational experiments were conducted on an Intel Core i7-12700H CPU @ 2.30GHz machine with 32GB RAM.

\subsection{Performance of the Iterative Alignment Method}

 In this subsection, we compare our IAM solution method with a benchmark that uses Gurobi. This benchmark optimizes maintenance and routing decisions, which requires evaluating a series of dynamic maintenance cost functions. Unfortunately, the dynamic maintenance cost function does not have a closed-form solution and requires evaluation of multiple integrations.
To address this, we develop a lower bounding envelope for the dynamic maintenance cost function based on a piece-wise linear approximation.

\textbf{Obtaining Lower Bounds:} 
To obtain lower bounds for the optimal objective function values, we create a lower bounding envelope for the dynamic maintenance cost function by using its structural properties. 
The dynamic maintenance cost function in Equation \eqref{eq:MaintenanceCost} is convex within its reasonable operational range \citep{yildirim2016sensor}. 
Hence, we find a piecewise linear function, which is a lower bound for our original dynamic maintenance cost function. 
Our approach is as follows. Let $\mathcal{K}$ be the number of linear functions we use to create the piecewise linear approximation. We first choose time point $h_{k}$ for $k\in \{1, 2, \dots, \mathcal{K}\}$.
Then, for each $h_{k}$, we choose a line that is tangent to the dynamic maintenance cost function at that point. 
Let $l_{k}$ be the $y$-intercept and $s_k$ be the slope of the $k^{th}$ line. Since these lines are tangent, they lie below the dynamic maintenance cost function. 
Hence, the following holds:
 \begin{equation}
f_{t_{o}}(t) \geq l_k + s_k t, \quad \forall k \in \{1, 2, \dots, \mathcal{K}\}. \label{eq:linearize}
\end{equation}

Therefore, by modifying Constraint \eqref{eq:maintCost} in our original mathematical model using inequalities in equation \eqref{eq:maintCost2}, we compute lower bounds for the incurred maintenance cost $\gamma$, which is exact at time $t=\pi$, and forms a valid lower bound for other time periods.
\begin{equation}
\gamma \geq  l_k + s_k \pi - My, \quad \forall k \in \{1, 2, \dots, \mathcal{K}\}. \label{eq:maintCost2}
\end{equation}
By solving the modified model using Gurobi, we obtain a lower bound for our original problem. Evidently, when we compare our solution results, we are, in fact, comparing our method with a lower bounding function of the original problem. Hence, the reported optimality gaps are upper bounds on the actual optimality gaps.

\textbf{Comparison of IAM with Gurobi-Based Lower Bounds:}
In the first part of our computational experiments, we use instances with up to 80 nodes since Gurobi fails to find any solutions beyond this limit. The names of the test instances are listed in Table \ref{table:performance_comparison} in the Dataset column. The column N/MN shows the number of nodes and the number of maintenance-capable nodes. For the computations, we enforce a two-hour time limit. 
In Table \ref{table:performance_comparison}, for IAM and Gurobi, we report the objective function values of the solutions found ($z$), lower bound ($LB$) and upper bound ($UB$), computation time in seconds ($t(sec)$), and the percentage gap between the lower and the upper bounds ($\%$ Gap). For instances where Gurobi cannot find the optimal solution, we report the value of the feasible solution.

\begin{table}[ht]
\centering
\small
\caption{Comparison of Iterated Alignment Method (IAM) with benchmark (Gurobi).}
\resizebox{\textwidth}{!}{%
\begin{tabular}{cc|ccccc|cccc|c}
\toprule
&  & \multicolumn{5}{c|}{\textbf{Gurobi}} & \multicolumn{4}{c|}{\textbf{IAM}} & \%Gap \\
\cmidrule(lr){3-7} \cmidrule(lr){8-11}
Dataset & N/MN & $z$ & Feasible & $LB$ & $t(sec)$ & \%Gap & $z$ & $UB$ & $LB$ & $t(sec)$ & \makecell{\footnotesize{$\frac{z_{\text{IAM}} - z_{\text{Gurobi}}}{z_{\text{Gurobi}}}$}} \\
\midrule
n10w180.002 & 10/2  & 2615.7 & - & - & 0.25 & 0\% & 2615.7 & 2622.3 & 2615.7 & 55.45 & 0.00\% \\
n10w10.002 & 10/3  & 2349.3 & - & - & 0.54 & 0\% & 2350.3 & 2350.3 & 2349.3 & 45.66 & 0.04\% \\
n20w180.003 & 20/3  & - & 2817.3 & 2572.5 & 7200 & 8.69\% & 2774.1 & 2775.9 & 2772.8 & 234.14 & - \\
n20w200.001 & 20/3  & - & 2241.6 & 2169.3 & 7200 & 3.22\% & 2227.0 & 2227.8 & 2216.6 & 200.93 & - \\
n40w120.001 & 40/4  & - & - & - & 7200 & - & 4378.4 & 4381.5 & 4378.4 & 991.49 & - \\
n40w120.002 & 40/4  & - & - & - & 7200 & - & 4395.9 & 4395.9 & 4392.8 & 407.31 & - \\
n80w100.001 & 80/6  & - & - & - & 7200 & - & 5317.5 & 5317.5 & 5294.1 & 4373.82 & - \\
\bottomrule
\end{tabular}
}
\label{table:performance_comparison}
\end{table}

\normalsize

Among these instances, Gurobi can find the optimal solution for instances up to 10 nodes. The IAM also finds the optimal solutions for these instances, although with slightly more computation time, but still within a minute. 
In datasets with 20 nodes, Gurobi fails to find exact solutions within the 7200-second limit and reports the current feasible solution and the lower bound.
In instances with 40 and 80 nodes, Gurobi cannot generate a feasible solution for the problem within the two-hour time limit.
However, in these instances, IAM finds optimal or near-optimal solutions, maintaining a 0.2\% gap between the upper and lower bounds. These results reveal that IAM significantly outperforms Gurobi in terms of the size of the instances that can be handled, and it finds solutions within reasonable time.

\begin{figure}[htbp]
    \centering
    \includegraphics[width=0.5\linewidth]{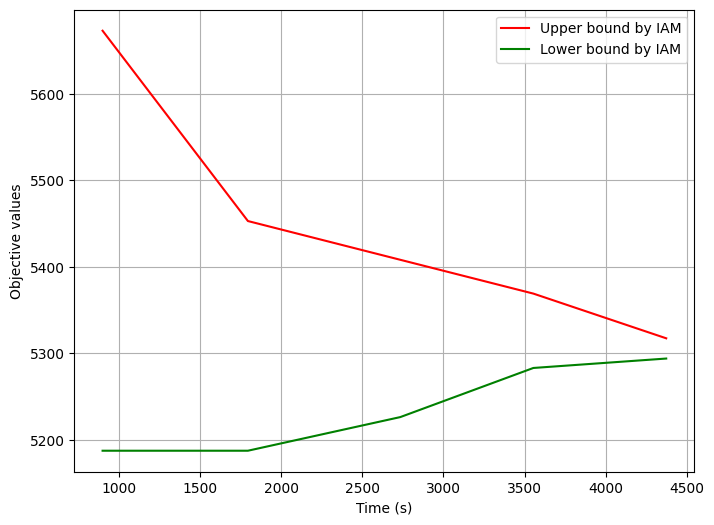} 
    \caption{Computational bounds in IAM for n80w100.001.}
    \label{fig:Comparison of Computational Bounds case:n80w100.001}
\end{figure}

The efficiency of IAM builds on the effectiveness of the lower and upper bounds that are computed through time-window subintervals for the maintenance-capable nodes. Using this approach, we can quickly reduce the feasible space and eliminate suboptimal solutions. This can be observed in Figure \ref{fig:Comparison of Computational Bounds case:n80w100.001}, where we display the upper and lower bounds computed by IAM for n80w100.001 instance with six maintenance nodes. We observe that our upper and lower bounds converge efficiently, while Gurobi fails to identify any bounds or feasible solutions within the observed time frame.

Our observations from the first part of the computational experiments underscore the potential of IAM to find high-quality solutions within reasonable time limits, especially for larger problem instances where exact solution methods become computationally infeasible. 
In the next part, we compare our sensor-driven approach to maintenance decisions with the traditional periodic review approach.

\subsection{Comparison of Sensor-Driven and Periodic Maintenance Approaches}
When integrating maintenance and routing decisions, we use a sensor-driven approach to model the relationship between these two decisions, where information on the current status of the vehicle obtained from the sensor data plays a key role. While this approach introduces higher computational challenges, it captures the dynamics better. An alternative approach is the traditional periodic maintenance, where maintenance decisions are made based on static measures, such as the age of the vehicle, ignoring its actual health status. In this section, we compare our sensor-driving routing and maintenance decisions (SDM-TSPTW) with those obtained under the periodic maintenance approach (PM-TSPTW) in terms of cost, the number of failures, and responsiveness to the flexibility (defined as the capability to adapt to different number of maintenance-capable nodes).
Different from the first part, our comparison of the sensor-driven and periodic maintenance approaches is based on actual operations. That is, once we determine the routing and maintenance decisions using both approaches, we evaluate them under the actual failure data obtained from \cite{gebraeel2005residual}. 
Before we present our computational results, we first describe how we compute the corresponding solutions using the periodic maintenance approach.

The periodic maintenance approach relies on pre-determined time windows, ignoring the health status of the equipment at the time. 
To determine routing and maintenance decisions using the periodic maintenance approach, we modify our IAM by using a new constraint to ensure that preventive maintenance occurs when the vehicle's age falls within predefined time periods.
According to the data in \cite{gebraeel2005residual}, the average time until failure is observed to be around 125 time units. Hence, we choose the periodic maintenance time window to be [100, 112], reflecting the 80\% to 90\% of average failure time.

Once we obtain solutions using sensor-driven and periodic maintenance approaches, we simulate these solutions based on vehicle degradation and failure data. 
If the vehicle does not break down before the planned maintenance, then the maintenance is performed as planned. This is called preventive maintenance. However, if the vehicle breaks down before the planned maintenance, then it needs to go through corrective maintenance, which is commonly much more expensive than preventive maintenance. In our computations, we use \$1000 and \$4000 for the preventive and corrective maintenance operations, respectively.

\textbf{Comparison of SDM-TSPTW and PM-TSPTW in terms of costs:}
Our first comparison of the SDM-TSPTW and PM-TSPTW is based on the actual cost realizations. In these experiments, we obtain solutions for different instances and simulate these solutions under 10 randomly selected scenarios. Then, we report the average total cost, routing cost, and maintenance cost for each instance.
We present our results in Table \ref{table:sensor_vs_periodic}.
Based on the results in Table \ref{table:sensor_vs_periodic}, we observe that the sensor-driven-based approach significantly outperforms the periodic maintenance-based approach in total cost, reducing the total cost of operation between 7.2\% and 27.1\%. The maintenance-related costs under the sensor-driven approach are 52\% less on average than the periodic maintenance approach, while the routing costs are similar, with some cases where the routing cost is larger. This is because both methods make integrated routing and maintenance decisions, the periodic maintenance approach, unlike the sensor-driven counterpart, does not consider the additional wear and tear of the routing decisions on the vehicle degradation level.

\begin{table}[h]
\centering
\small
\caption{Comparison of SDM-TSPTW and PM-TSPTW in terms of costs.}
\resizebox{\textwidth}{!}{%
\begin{tabular}{cc|ccc|ccc|c}
\toprule
 &  & \multicolumn{3}{c|}{\textbf{SDM-TSPTW}} & \multicolumn{3}{c|}{\textbf{PM-TSPTW}} & \makecell{\% Cost \\ Reduction } \\
\cmidrule(lr){3-5} \cmidrule(lr){6-8}
Dataset & N/MN & \makecell{Total\\ Cost} & \makecell{Routing \\ Cost} & \makecell{Maintenance \\ Cost} &  \makecell{Total \\Cost} & \makecell{Routing \\ Cost} & \makecell{Maintenance \\ Cost} &  $\frac{PM}{SDM}-1$ \\
\midrule
n10w180.002 & 10/2 & 3805.6 & 2505.6 & 1300.0 &  4405.6 & 2505.6 & 1900.0 &  15.8 \\
n10w10.002 & 10/3 & 3194.8 & 2044.8 & 1150.0 &  4059.2 & 2059.2 & 2000.0 &  27.1 \\
n20w180.003 & 20/3 & 3869.2 & 2419.2 & 1450.0 &  4146.0 & 2196.0 & 1950.0 &  7.2 \\
n20w200.001 & 20/3 & 3107.9 & 1907.9 & 1200.0 &  3907.9 & 1907.9 & 2000.0 &  25.7 \\
n40w120.001 & 40/4 & 5203.5 & 4053.5 & 1150.0 &  5653.5 & 4053.5 & 1600.0 &  8.6 \\
n40w120.002 & 40/4 & 5518.0 & 4068.0 & 1450.0 &  6368.0 & 4068.0 & 2300.0 &  15.4 \\
\bottomrule
\end{tabular}
}
\label{table:sensor_vs_periodic}
\end{table}

\textbf{Comparison of SDM-TSPTW and PM-TSPTW in terms of number of failures:}
In this part, we turn our attention to the reliability aspect of the two maintenance strategies by examining the incidence of failures that occur under sensor-driven and periodic maintenance strategies. 
To investigate this, we conduct a series of 100 simulations for each of the instances considered in our study.
In each simulation run, we use randomly selected sensor data based on which we derive a dynamic maintenance cost and make decisions for the SDM-TSPTW. We also make the decisions for the PM-TSPTW as described earlier. 
While operating with these solutions, if the scheduled maintenance time 
exceeds the actual failure time of the vehicle, it results in vehicle failure, requiring corrective maintenance. 
The results of these simulations are in Table \ref{tab:failures}, where we present the number of failures that occur under both maintenance approaches.  
We observe that SDM-TSPTW results in a much lower number of failures compared to the PM-TSPTW, with an average of 9.17 failures per instance compared to that of 31.5 for the PM-TSPTW in 100 runs.
These results reinforce that the sensor-driven strategy is more effective in creating more reliable operations planning, which not only reduces the cost but also the secondary challenges associated with vehicle breakdowns, such as customer dissatisfaction.

\begin{table}[h]
    \centering
    \caption{Comparison of SDM-TSPTW and PM-TSPTW in terms of number of failures.}
    \label{tab:failures}
    \small
    \begin{tabular}{cc|c|cc}
        \hline
        \makecell{\\Dataset} & \makecell{\\N/MN} & \multicolumn{1}{c|}{\textbf{SDM-TSPTW}} & \multicolumn{1}{c}{\textbf{PM-TSPTW}} \\
        \cmidrule(lr){3-3} \cmidrule(lr){4-4}
         &  & \# of failures & \# of failures \\
        \hline
        n10w180.002 & 10/2 & 10 & 32 \\
        n10w10.002 & 10/3 & 5 & 33 \\
        n20w180.003 & 20/3 & 15 & 32 \\
        n20w180.004 & 20/3 & 5 & 33 \\
        n40w120.001 & 40/4 & 5 & 22 \\
        n40w120.002 & 40/4 & 15 & 42 \\
        \hline
    \end{tabular}
\end{table}

\textbf{Comparison of SDM-TSPTW and PM-TSPTW under different maintenance flexibility:}
Finally, we examine how the sensor-driven and periodic maintenance strategies respond to increased flexibility in the system, where we associate flexibility with the number of maintenance-capable nodes. 

To investigate this, create instances with different numbers of maintenance-capable nodes by modifying the instance n40w120.001. Among 40 nodes, we initially choose one of them as the maintenance-capable node by solving the $p$-median problem. Then, we add more maintenance-capable nodes to this set, up to 7 nodes. After obtaining solutions using the sensor-driven and periodic maintenance approaches, we simulate the results on 10 different sensor data realizations. In Table \ref{tab:flexibility}, we present the average total cost, routing cost, and maintenance cost over 10 realizations for each number of maintenance nodes.

\begin{table}[h]
\centering
\small
\caption{Comparison of SDM-TSPTW and PM-TSPTW under different numbers of maintenance-capable nodes. 
The costs are the averages of 10 random realizations.}
\resizebox{\textwidth}{!}{%
\begin{tabular}{c|ccc|ccccccc}
\toprule
\# of Maintenance Nodes & \multicolumn{3}{c|}{\textbf{SDM-TSPTW}} & \multicolumn{3}{c}{\textbf{PM-TSPTW}} \\
\cmidrule(lr){2-4} \cmidrule(lr){5-7}
&  Total Cost & \makecell{Routing \\ Cost} & \makecell{Maintenance \\ Cost}  & Total Cost & \makecell{Routing \\ Cost} & \makecell{Maintenance \\ Cost}  \\
\midrule
1 & 5353.6 & 4053.6 &  1300.0 & 7153.6 & 4053.6 & 3100.0  \\
2 & 5353.6 & 4053.6 &  1300.0 & 5635.6 & 4053.6 & 1900.0  \\
3 & 5053.6 & 4053.6 &  1000.0 & 5953.6 & 4053.6 & 1900.0  \\
5 &  4918.9 & 3908.9 &  1000.0 & 5475.5 & 3575.5 & 1900.0  \\
7 &  4874.0 & 3874.0 &  1000.0 & 4875.5 & 3575.5 & 1300.0 \\
\bottomrule
\end{tabular}
}
\label{tab:flexibility}
\end{table}

In real-world applications, transforming a node into a maintenance-capable node may require considerable investments, including setting up facilities, training personnel, acquiring tools and equipment, and integrating logistics and information technology systems for efficient operations. These investments are substantial and must be justified by a corresponding reduction in maintenance and operational costs. 
Thus, we conclude that the SDM-TSPTW not only ensures greater operational reliability, as evidenced by the significantly smaller number of failures in our simulations, but it also represents a more economically prudent option. Sensor-driven approach to making maintenance decisions for vehicles performing deliveries efficiently predicts the vehicle breakdown and allows us to compute the cost of maintenance considering the tradeoff between the early and late maintenance operations. Hence, it also reduces the need for extensive maintenance infrastructure.


\section{Summary}
\label{summary}

In this paper, we consider simultaneous single-vehicle routing with time windows and maintenance decisions and develop a sensor-driven predictive maintenance-based approach to capture the impact of routing decisions on vehicle degradation levels. 
We call this problem the Sensor-Driven Traveling Salesman Problem with Time Windows (SDM-TSPTW). Our modeling approach is based on integrating dynamic maintenance cost function that is derived based on vehicle sensor data.
Our proposed mixed integer programming (MIP) model addresses operational efficiency by optimizing routing decisions but also significantly contributes to the longevity and reliability of the vehicle through timely maintenance interventions, considering the tradeoffs between early and late maintenances.
 However, due to the nature of the dynamic maintenance cost function, solving this problem using commercial solvers is intractable. Hence, we have develop a new solution method, called the Iterative Alignment Method (IAM), to solve the problem. 
We compare our proposed SDM-TSPTW approach with benchmarks to measure its performance in different settings. First, we show that IAM outperforms Gurobi in terms of both solution quality (Gurobi cannot find solutions within the time limit) and computation time. 
In addition, we compare our sensor-driven approach to maintenance decisions with the periodic maintenance approach, which ignores the current status of the vehicle and uses a static metric to decide when to perform maintenance. We show that our sensor-driven approach is significantly better than the periodic maintenance-based approach, in terms of the total cost, number of vehicle failures, and leveraging the flexibility in the systems in terms of the number of maintenance-capable locations.

Our work also opens up future research directions. Studying routing problems with multiple vehicles is an important direction. Especially when we can reassign the jobs of a vehicle that breaks down during the route. In our study, we assume that the route is fully determined before dispatching, and there is no rerouting. In settings where the route decisions are made dynamically, the sensor data can be used more efficiently.

\section{Acknowledgments}
The authors are thankful to the National Science Foundation for supporting this work through grant \#2104455.

\appendix
\vspace{-0.3cm}

\section*{Appendix A: Proof of Propositions}
    \vspace{-0.3cm}
\label{appendix:proof}

\begin{proof}[Proof of Proposition 1]
The proof is straightforward. 
For contradiction, assume that the optimal maintenance solution resides within subinterval \(q\) of node \(i\). Since \(C^r \tau^i_q + \underline{g}^i_q > C^r \tau^i_k + \overline{g}^i_k\), the lowest possible cost in subinterval \(q\) (\(C^r \tau^i_q + \underline{g}^i_q\)) is greater than the highest possible  cost in subinterval \(k\) (\(C^r \tau^i_k + \overline{g}^i_k\)). 
Therefore, choosing subinterval \(q\) cannot be optimal as it leads to a higher cost than that of subinterval $k$. 
\end{proof}

\begin{proof}[Proof of Proposition 2]
Consider the upper bound of the maintenance cost in subinterval $q$ of maintenance node $i$, $\bar{g}^i_q$.
If the maintenance is performed in subinterval \(q\), the upper bound on the total cost is \(C^r \tau^i_q + \overline{g}^i_q \).  
Since \(C^r \tau^i_q + \overline{g}^i_q \) is the upper bound for all feasible solutions when the maintenance is performed at node $i$ in time window interval $q$, an upper bound for the SDM-TSPTW when the maintenance is performed at node $i$ is \(\min_{q \in S_i} \{C^r \tau^i_q + \overline{g}^i_q \}\).
\end{proof}

\begin{proof}[Proof of Proposition 3]
Lower bound of the optimal maintenance cost in subinterval $q$ of maintenance node $i$ is $\underline{g}^i_q$.
If we consider a feasible solution where the maintenance is performed in subinterval $q$ of node $i$, then the corresponding objective function value is greater than or equal to \(C^r \tau^i_q + \underline{g}^i_q\).
Hence, considering all subintervals of node $i$, \(\min_{q \in S_i}\{C^r \tau^i_q + \underline{g}^i_q \})\) is the minimum cost that can be achieved if the maintenance occurs at node $i$.
\end{proof}

\begin{proof}[Proof of Proposition 4]

The proof has two parts. The first part proves that the solution converges to an $\epsilon$ optimal solution. The second part shows finite convergence. Before we present the proof, we will first provide new notation for ease of explanation.

Let us assume that the optimal solution for the problem has a maintenance at the $i^{th}$ maintenance capable node at time $t^*$. 
At the $v^{th}$ iteration, assume the time $t^*$ lies within period $q^*$.
For ease of notation, let $\tau^{i}_{t}$ be the corresponding optimal route duration when the maintenance is performed at node $i$ at time period $t$.
In addition, let us define the following: (i) upper bound for optimal solution $\bar{o} = C^r \tau_{q^*}^i + \bar{g}^i_{q^*}$, (ii) actual optimal solution $o = C^r \tau^{i}_{t^*} + g^i_{t^*} \leq \bar{o}$, and (iii) lower bound for optimal solution $\underline{o} = C^r \tau_{q^*}^i + \underline{g}^i_{q^*}$. Further, let us also denote a set of feasible solutions that are not $\epsilon$-optimal (with a total cost higher than the optimal total cost plus $\epsilon$) as $\mathcal{F}$.
For every element $e$ of this set $\mathcal{F}$, we denote the corresponding maintenance capable node as $j_e$, and time period as $k_e$, and define the following: (i) upper bound for feasible solution $e$ is $\bar{f}_e = C^r \tau_{k_e}^{j_e} + \bar{g}^{j_e}_{k_e}$, (ii) actual cost for the feasible solution $e$ is $f_e = C^r \tau^{j_e}_{t_e} + g^{j_e}_{t_e}$, and (iii) lower bound for feasible solution $e$ is $\underline{f}_e = C^r \tau_{k_e}^{j_e} + \underline{g}^{j_e}_{k_e}$. From the definition of $\mathcal{F}$, the following always holds: $f_e \geq o+\epsilon, \forall e \in \mathcal{F}$.

\underline{\textit{Proof for $\epsilon$-optimality:}} Let us assume for contradiction that the algorithm stops (i.e. $U-L \leq \epsilon + 2 \delta^v$), and one of the feasible solutions that are not $\epsilon$ optimal, $e\in \mathcal{F}$ is not dominated. For any feasible solution $e$, we can make the following assertions: For the upper bound $U$, we can claim: $ f-\delta^v \leq \underline{f}_e \leq U$
, because if the upper bound for any solution is less than $\underline{f}_e$, then the solution $e$ would be eliminated by Proposition 1. Lower bound for the problem is the minimum of the lower bound for every non-dominated solution. The smallest actual solution is the optimal solution ${o}$, hence the corresponding lower bound has the following property $
L \leq \bar{o} \leq o+\delta^v$. Evidently, if $\epsilon + 2\delta^v \geq U-L$ then $\epsilon + 2\delta^v \geq U-L \geq (f-\delta^v) - (o + \delta^v) = f-o - 2\delta^v$, hence $f - o \leq \epsilon$. However, our initial assumption was $f - o \geq \epsilon$, which proves that $f$ cannot be a part of the non-dominated solution when the algorithm converges. Since every feasible solution within $\mathcal{F}$ is dominated at convergence, the best solution is $\epsilon$-optimal.

\underline{\textit{Proof for finite convergence:}} Each iteration $v$, has the corresponding $\delta^v = \delta^o / (b^v)$. For every solution $f_e$ that is not $\epsilon$-optimal defined by the set $\mathcal{F}$, we define a term $\gamma_e>0$, such that $f_e = o + \epsilon + \gamma_e$. Evidently, $\underline{f}_e \geq f_e - \delta^v = o + \epsilon + \gamma_e - \delta^v$. For any $\gamma_e >0$, the sufficient condition $f_e - \delta^v = o +\epsilon + \gamma_e - \delta^v \geq o + \delta^v$ implies that $\underline{f}_e > \bar{o}$. This follows from the fact that in the sufficient condition the first and second term is replaced by the corresponding lower and upper bound, respectively. This sufficient condition can be satisfied if $\epsilon + \gamma_e  \geq 2 \delta^v$, and hence $\delta^v \leq (\epsilon + \gamma_e)/2$, which is guaranteed to be achieved at iteration $v \geq log(2\delta^0/(\epsilon+\gamma_e))/log(b)$ for any solution that is not $\epsilon$-optimal. The maximum of these finite iteration limits for different solutions yields the number of iterations in the algorithms, which would also be finite.
\end{proof}

\hspace{-6mm}\textbf{\textit{Data Availability Statement:}} Authors have permission to use sensor data in numerical study, but cannot make it publicly available. TSPTW instance data is available online.


{
\bibliographystyle{apalike}
\newlength{\bibitemsep}\setlength{\bibitemsep}{0\baselineskip plus .05\baselineskip minus .05\baselineskip}
\newlength{\bibparskip}\setlength{\bibparskip}{0pt}
\let\oldthebibliography\thebibliography
\renewcommand\thebibliography[1]{%
  \oldthebibliography{#1}%
  \setlength{\parskip}{\bibitemsep}%
  \setlength{\itemsep}{\bibparskip}%
}
\bibliography{Arxiv2}

\end{document}